\documentclass[12pt, a4paper, reqno]{amsart}
\usepackage{amsmath, amsthm, amsfonts, amssymb, xcolor}
\usepackage{tikz}
\usetikzlibrary{arrows.meta, positioning, cd}
\usepackage[hidelinks]{hyperref}
\usepackage[margin=1in]{geometry}
\usepackage{enumitem}

\newtheorem{theorem}{Theorem}[section]
\newtheorem{lemma}[theorem]{Lemma}
\newtheorem{proposition}[theorem]{Proposition}
\newtheorem{corollary}[theorem]{Corollary}
\newtheorem{fact}[theorem]{Fact}
\theoremstyle{definition}
\newtheorem{definition}[theorem]{Definition}
\newtheorem{example}[theorem]{Example}

\newtheorem{remark}[theorem]{Remark}
\numberwithin{equation}{section}

\newenvironment{claim}{%
  \vspace{0.5em}\hfill \break\noindent \textbf{Claim.}
}

\newenvironment{claimproof}{%
  \par\noindent \textit{Proof: }%
}{\space\hfill $\qed_{\operatorname{Claim}}$\vspace{0.5em}\newline}

\newcommand{\gtcvf}{G\operatorname{TCVF}}
\newcommand{\gtcf}{G\operatorname{TCF}}
\newcommand{\cO}{\mathcal{O}}
\newcommand{\cL}{\mathcal{L}}
\newcommand{\s}{\sigma}
\newcommand{\ec}{\prec_\exists}
\newcommand{\ecoag}{\ec^{\operatorname{oag}}}
\newcommand{\ecval}{\ec^v}
\newcommand{\ecg}{\ec^G}
\newcommand{\ecgval}{\ec^{v,G}}
\newcommand{\ecrng}{\ec^{\operatorname{rng}}}
\newcommand{\Kv}{\left( K, v\right)}

\newcommand{\Z}{\mathbb{Z}}

\newcommand{\Kvs}{\left( K, v,\s\right)}
\newcommand{\lng}{\cL_{\cO, G}}
\newcommand{\lngg}{\cL_{G}}
\newcommand{\lngval}{\cL_{\cO}}
\newcommand{\lngring}{\cL_{\operatorname{ring}}}
\newcommand{\ntp}{\operatorname{NTP}_2}
\newcommand{\bdn}{\operatorname{bdn}}

\begin{document}

\title{A note on valued fields with finite group actions}

\author[P. Błaszkiewicz, J. Gogolok]{Piotr Błaszkiewicz$^1$ and Jakub Gogolok$^2$}

\thanks{$^{1}$ Supported by the Narodowe Centrum Nauki grant no. 2021/43/B/ST1/00405.}
\thanks{
$^{2}$ Supported by the Narodowe Centrum Nauki grant no. 2023/49/N/ST1/02512.
}

\address{$^{1}$ Instytut Matematyczny, Uniwersytet Wrocławski, Wrocław, Poland}
\email{piotr.blaszkiewicz@math.uni.wroc.pl}

\address{$^{2}$ Instytut Matematyczny, Uniwersytet Wrocławski, Wrocław, Poland}
\email{jakub.gogolok@math.uni.wroc.pl}
\urladdr{http://www.math.uni.wroc.pl/\textasciitilde gogolok/}



\begin{abstract}
We consider valued fields equipped with an action of a given finite group by automorphisms. We show that the theory of such structures admits a model companion which is axiomatized in an Ax-Kochen/Ershov fashion. We prove this model companion has $\ntp$ and that the valuation is definable in the pure field language.
\end{abstract} \maketitle

\tableofcontents

\section{Introduction}
\noindent\textbf{Summary of main results.}
In this note we prove that the theory of valued fields equipped with an action of a given finite group $G$ (we call such structures \textbf{valued $G$-fields}, see Definition \ref{def: gvf}) has a model companion. Our main result (Theorem \ref{thm: main ake}) states that existentially closed valued $G$-fields admit a characterization in the spirit of Ax-Kochen-Ershov. From this the existence of a model companion follows (see Theorem \ref{thm: mc}). We denote this model companion by $\gtcvf$. We show that the valuation in models of $\gtcvf$ is definable in the language of rings (see Proposition \ref{prop: def val}) and that certain extensions of $\gtcvf$ have $\ntp$ (see Theorem \ref{thm: ntp2}). We also produce some interesting examples of model complete pure fields (see Proposition \ref{prop: mc field} and Remark \ref{rem: mc fields}).

\vspace{0.5em}
\noindent\textbf{Motivations.} A prototypical example of a model companion of a theory of fields is the theory of algebraically closed fields ACF, being the model companion of the theory of (pure) fields.\footnote{To speak about a model companion one has to fix a language in which we are considering structures. We usually skip the language where there is a natural ``default language'' know from the context. See also Subsection \ref{subsec: langs} for a broader discussion of this.} The existence of a model companion $T^*$ of a theory $T$ implies the existence of a model $\mathcal{U}\models T$ which is both existentially closed and sufficiently saturated. Such a model generalizes the notion of a ``universal domain'' appearing in Weil's foundations of algebraic geometry. Thus, the existence of $T^*$ yields a framework for ``$T$-algebraic geometry'' and this is especially interesting then $T$ is a theory of fields (with possibly some additional structure). Prominent examples of $T^*$ include the theories of
\begin{enumerate}
\item differentially closed fields DCF (see \cite{DCF0}, \cite{DCFp}),
\item existentially closed difference fields ACFA (see \cite{ACFA}),
\item algebraically closed valued fields ACVF (see \cite{ACVF}).
\end{enumerate}
The theory ACFA can be viewed as the model companion of the theory of fields with a $\mathbb{Z}$-action. In \cite{GTCF} Hoffmann and Kowalski introduced a theory $\gtcf$ and showed that it is the model companion of the theory of fields with an action of a given finite group $G$ (see Fact \ref{gtcf} and Fact \ref{gtcf2}). Chernikov and Hills introduced in \cite{VFA} the model companion of the theory of (contractive) difference valued fields VFA, so in other words: the model companion of the theory of valued fields with (contractive) $\mathbb{Z}$-action. The relationship betweens the theories we consider and the theories mentioned above can be summarized in the following diagram.
\begin{center}
\begin{tikzpicture}
\node[draw] (acf) at (-4.5,0) {ACF};
\node[draw] (acfa) at (-6.5,-2) {ACFA};
\node[draw] (gtcf) at (-2.5,-2) {$\gtcf$};

\draw[->, thick] (acf) -- (acfa);
\draw[->, thick] (acf) -- (gtcf);

\node[draw] (acfv) at (4.5,0) {ACFV};
\node[draw] (gtcvf) at (6.5,-2) {$\gtcvf$};
\node[draw] (vfa) at (2.5,-2) {VFA};

\draw[->, thick] (acfv) -- (vfa);
\draw[->, thick] (acfv) -- (gtcvf);

\draw[->, thick] (-1.8,-1) -- (2,-1);

\end{tikzpicture}
\end{center}

\vspace{0.5em}

\noindent\textbf{Conventions.} Throughout this note we work inside a big algebraically closed field $\Omega$ equipped with a nontrivial valuation $v$. We also fix a \textbf{nontrivial} and \textbf{finite} group $G$. Some of our results are also valid for infinite $G$ but most use finiteness in an essential way.
\vspace{0.5em}

\noindent\textbf{Acknowledgment.} We would like to thank Piotr Kowalski and Franz-Viktor Kuhlmann for numerous valuable discussions. We also thank the Wrocław model theory group for various comments given during the seminar talks given by the first author at the University of Wrocław.

\section{Preliminaries}

In this section, we first set the languages we will work with. We also discuss some standard facts about fields with finite group actions and valued fields. Finally, we prove a few lemmas about these objects.

\subsection{Of languages and reducts}\label{subsec: langs}
We take now the opportunity to fix some notation and clarify certain model-theoretic nuances appearing in our work.

We are interested in three types of structures
\begin{enumerate}
    \item $G$-fields (fields with an action of a finite group $G$, see Definition \ref{def: gfield}), 
    \item valued fields,
    \item valued $G$-field (a $G$-field with a valuation preserved by the group action, see Definition \ref{def: gvf}).
\end{enumerate}
Unless necessary, we skip the additional structure on a field in notation and write, e.g. ``$K$ is a valued field'' instead of ``$(K,v)$ is a valued field''. In such a case, the valuation is always denoted by $v$, the action of an element $g\in G$ is denoted by $\s_g$ and tuple $\left\langle \s_g\colon g\in G\right\rangle$ is denoted by $\s$.

We are mostly interested in existentially closed (sub)structures of the above types. The phrase (say) ``$K$ is existentially closed as a valued $G$-field'' is meaningless unless we specify the language over which we consider a valued $G$-field. In Table 1. we specify the languages for all the structures we consider in this paper. We also give a notation for the relation ``$M$ is existentially closed (e.c.) in $N$''.
\begin{table}[h!]
\label{table}
\begin{tabular}{|l|l|c|}
\hline
\multicolumn{1}{|c|}{\textbf{Structure}} & \multicolumn{1}{c|}{\textbf{Language}}                   & \textbf{``$M$ is e.c. in $N$''} \\ \hline
ordered abelian groups                   & $\left\{ +,-,0, \le\right\}$                             & $M \ecoag N$                    \\ \hline
rings/fields                             & $\lngring=\left\{ +,-,\cdot, 0, 1\right\}$               & $M \ecrng N$                    \\ \hline
valued fields                            & $\lngval=\lngring\cup\left\{ \cO \right\}$               & $M \ecval N$                    \\ \hline
$G$-fields                               & $\cL_{G}=\lngring\cup\left\{ \s_g\colon g\in G \right\}$ & $M \ecg N$                      \\ \hline
valued $G$-fields                        & $\lng=\lngval \cup \cL_{G}$                              & $M \ecgval N$                   \\ \hline
\end{tabular}
\vspace{0.5em}
\caption{Specification of languages and notation}
\end{table}

The language $\lngval$ is one-sorted and $\cO$ is a unary predicate for the valuation ring. Instead of $\lngval$ we could have chosen e.g. the standard three-sorted language $\cL_{\operatorname{val}}$. Of course those two languages (with their intended interpretations) yield the same interpretable sets in valued fields and the same notions of substructures, but a priori may lead to different notions of existentially closed extensions. The following fact says that it is not the case. We leave its proof to the Reader.
\begin{lemma}
Suppose that $K\subseteq L$ is an extension of valued fields. The following are equivalent:
\begin{enumerate}
    \item $K$ is e.c. in $L$ as an $\lngval$-structure,
    \item $K$ is e.c. in $L$ as an $\cL_{\operatorname{val}}$-structure.
\end{enumerate}
\end{lemma}

\subsection{\texorpdfstring{$G$}{G}-fields}
In this section, we wish to define what a $G$-fields are, as well as introduce notions and results we will use in the further parts of this paper. The following are mostly based on \cite{GTCF}. 

\begin{definition}\label{def: gfield}
A \textbf{$G$-field} is a field $K$ equipped with an action of $G$ by field automorphisms. A $G$-field is called \textbf{strict} if the action of $G$ is faithful. We call the fixed field $C_K:=\left\{x\in K\colon(\forall g\in G)(gx = x)\right\}$ the \textbf{field of constants} of $K$.
\end{definition}
The following easy fact will be used throughout this note.
\begin{lemma}\label{lemma: G compositum}
Let $K$ be a $G$-field and let $C'$ be an extension of $C_K$ linearly disjoint from $K$ over $C_K$. Then the action of $G$ has a unique extension to $KC'$ for which the field of constants is $C'$.
\end{lemma}
\begin{proof}
Uniqueness is evident, so let us prove existence. By the linear disjointness and the finiteness of $G$ we have that the compositum $KC'$ is isomorphic (as a $C_K$-algebra) to the tensor product $K\otimes_{C_K} C'$. The group $G$ naturally acts on this tensor product and hence on $KC'$; explicitly, we have
$$\widetilde{\s}_g\left( \sum_{i=1}^n x_i\otimes c_i \right) := \sum_{i=1}^n \s_g\left( x_i\right)\otimes c_i,$$
for $g\in G$, $x_1,\dotsc, x_n\in K$, $c_1,\dotsc,c_n\in C'$. It is easy to check that the field of constants of this action is precisely $C'$.
\end{proof}

\begin{lemma}\label{lemma: perfect constants}
Let $L$ be a $G$-field of characteristic $p>0$ and assume that every element of $C_L$ has a $p$th root in $L$. Then $C_L$ and $L$ are perfect.   
\end{lemma}
\begin{proof} Take $a\in C_L$ and $b\in L$ such that $b^p = a$. Then, for any $g\in G$ we have that 
$$\sigma_g\left( b \right) ^p = \sigma_g\left( b^p \right) = \sigma_g\left( a \right) = a = b^p,$$
thus $\sigma_g\left( b \right) = b$. Therefore $b\in C_L$ and hence $C_L$ is perfect. Since the extension $C_L\subseteq L$ is algebraic, the field $L$ is perfect, too.
\end{proof}

\begin{remark}\label{rem: K definable in C}
Let $K$ be a $G$-field. We claim that the $\lngg$-structure $K$ is interpretable in the $\lngring$-structure $C_K$. To see this note that $K\cong C_K^{\left|G\right|}$ as $C$-vector spaces and under this isomorphism the multiplication on $K$ is encoded by a set of $\left| G \right|$-matrices with coefficients in $C_K$ (see also \cite[Remark 2.3]{GTCF}). Note that this interpretation is given by quantifier-free formulas and no quotients are involved.
\end{remark}
Henceforth \textbf{we will always treat $K$ as a definable object in the structure $C_K$.} This allows us say e.g. that a subset $X\subseteq K$ is definable in $C_K$. 

An easy but very useful application of Remark \ref{rem: K definable in C} is the following lemma, whose proof we omit (see also \cite[Remark 2.3]{GTCF}).

\begin{lemma}\label{lemma: ec red to c pure}
Let $K\subseteq L$ be an extension of strict $G$-fields. Then $K \ecg L$ if and only if $C_K\ecrng C_L$.
\end{lemma}

Finally, we need the following results about existentially closed $G$-fields.
\begin{fact}[Theorem 3.29 in \cite{GTCF}]\label{gtcf}
    Let $K$ be a $G$-field. Then $K$ is existentially closed as a $G$-field if and only if $K$ is strict, $C_K$ is PAC and $K$ does not admit any proper algebraic $G$-field extensions. 
\end{fact}
\begin{fact}[Theorem 3.29 in \cite{GTCF}]\label{gtcf2}
The theory of $G$-fields admits a model companion $\gtcf$.    
\end{fact}

\subsection{Valued fields}

Let $\Kv$ be a valued field. We denote the value group of $\Kv$ by $vK$, the residue field by $Kv$ and the valuation ring by $\mathcal{O}_K$ (following the convention from \cite{FVKbook}). We denote by $\operatorname{res}\colon \cO_K\to Kv$ the residue map.

Recall that for a finite extension $C\subseteq K$ of valued fields we have the so-called \textbf{fundamental inequality} (see \cite[Theorem 3.3.4]{engler}):
\[
\left[ K : C \right] \ge \left[ vK: vC \right] \cdot [ Kv: Cv ].
\]

A finite extension $C\subseteq K$ of valued fields is called \textbf{defectless} if the above inequality is in fact an equality. On the other extreme are the \textbf{immediate extensions}, that is extensions $C\subseteq C'$ for which $vC=vC'$ and $Cv=C'v$. A valued field $(C,v)$ is called \textbf{algebraically maximal} if it has no proper algebraic immediate extension. A valued field $(C,v)$ is called \textbf{henselian}\ if the valuation $v$ extends uniquely to any algebraic extension of $C$.\footnote{For more equivalent definitions of henselianity we refer to Theorems 4.1.3 and 4.1.7 from \cite{engler}.} Every valued field $C$ admits an immediate extension $C\subseteq C^h$ which is henselian. This extension is unique up to isomorphism over $C$ and is known as the \textbf{henselianization} of $C$ (see \cite[Theorem 5.2.2, Theorem 5.2.5]{engler}). In particular, algebraically maximal valued fields are henselian. 
\begin{fact}[Lemma 2.4 in \cite{Kuhlmann}]\label{fact: immediate and defectless are ld}
Let $C\subseteq C'$ be an immediate
extension of valued fields. If $C\subseteq K$ is a finite defectless extension admitting a unique extension of the valuation, then the same holds for $C'\subseteq KC'$ and $K$ is linearly disjoint from $C'$ over $C$.
\end{fact}

\begin{lemma}[Lemma 11.81 in \cite{FVKbook}]\label{lemma: val indep are lin indep}
Let $C\subseteq C'$ be a valued field extension. If
\begin{enumerate}
\item $C\subseteq K$ is a finite defectless valued field extension,
\item $Kv$ and $C'v$ are linearly disjoint over $Cv$, and
\item $vK\cap vC'=vC$
\end{enumerate}
then $K$ is linearly disjoint from $C'$ over $C$ and the valuation $v$ extends uniquely from $C'$ to a valuation $\widetilde{v}$ on $KC'$. Furthermore the extension $(C',v)\subseteq (KC',\widetilde{v})$ is defectless, $\widetilde{v}(KC')=vK+vC'$ and $Kv.C'v=(KC')\widetilde{v}$.
\end{lemma}

For our purposes the notion of tame valued fields, introduced by Kuhlmann in \cite{tame}, will be crucial. 
\begin{definition}
A valued field $\Kv$ is called \textbf{tame} if it is algebraically maximal, $Kv$ is perfect and $vK$ is $p$-divisible, where $p$ is the characteristic of $Kv$ (for $p=0$ this means that $vK$ is divisible).
\end{definition}
The above is not the original definition of tame fields, but rather an equivalent one, which is part of the content of \cite[Theorem 3.2]{tame}. Tame fields are important because they form in a sense the largest class of valued fields for which one can hope for an Ax-Kochen/Ershov principle.
\begin{fact}[Theorem 1.4 in \cite{tame}]\label{fact: ake for tame}
Assume that $\Kv$ is a tame field and $K \subseteq L$ is an extension of valued fields. If $Kv\ecrng Lv$ and $vK\ecoag vL$, then $K\ecval L$.
\end{fact}

\section{Valued \texorpdfstring{$G$}{G}-fields}
The goal of this sections is to develop some basic algebraic results for the valued fields equipped with a finite group actions.
\subsection{Generalities}
This section introduces valued $G$-fields, that is, fields equipped with a valuation and an action $G$ by automorphisms preserving the valuation ring set-wise. We are interested in the interaction between these two structures, in particular in the induced action of $G$ on the value group and the residue field, as well as in the valuation on the field of constants.

The diagram below presents most reducts and auxiliary structures we will talk about while working with valued $G$-fields.
\begin{center}
\begin{tikzpicture}[
    >=Latex,
    node distance=0.7cm and 1.5cm,
    every node/.style={font=\normalsize}
]
% top
\node (A) {$\bigl(K,v,\sigma\bigr)$};

% first level
\node (C) [below=of A] {$\bigl(C_K,v\bigr)$};
\node (B) [left=of C] {$\bigl(Kv,\sigma\bigr)$};
\node (D) [right=of C] {$vK$};

% side nodes
\node (Z) [right=of A] {$\bigl(K,\sigma\bigr)$};
\node (T) [left=of A] {$\bigl(K,v\bigr)$};

% second level
\node (E) [below=of B] {$C_{Kv}$};
\node (F) [below=of C] {$C_Kv$};
\node (G) [below=0.85cm of D] {$vC_K$};

% arrows
\draw[->] (A) -- (B);
\draw[->] (A) -- (C);
\draw[->] (A) -- (D);
\draw[->] (A) -- (Z);
\draw[->] (A) -- (T);

\draw[->] (B) -- (E);
\draw[->] (C) -- (F);

\draw[->] (D) -- (G);
\draw[->] (C) -- (G);

\end{tikzpicture}
\end{center}
One may wonder why there is no $( vK,\s )$ in the diagram. This is explained by the following remark.
\begin{remark}
Let $\Kvs$ be a valued $G$-field. The action of $G$ on $K$ induces an action on the residue field $Kv$ by field automorphisms and on the value group $vK$ by order-preserving group automorphisms. The action of $G$ on $vK$ is by force trivial, due to the finiteness of $G$. Indeed, assume $g\in G$ has order $n$ and let $\gamma\in vK$ be such that $\s_g\left( \gamma \right) \neq \gamma$, say $\s_g\left( \gamma \right) > \gamma$. Since $\s_g$ is order-preserving, we have that
$$\gamma < \s_g\left( \gamma \right) < \s_g^2\left( \gamma \right)<\dotsc < \s_g^n\left( \gamma \right)=\gamma,$$
which is a contradiction.
\end{remark}

\begin{definition}\label{def: gvf}
A \textbf{valued $G$-field} is a valued field $\left( K, v \right)$ together with an action of $G$ on $K$ by valued field automorphisms, i. e. field automorphisms preserving the valuation ring $\mathcal{O}_K$ set-wise. A valued $G$-field is called \textbf{residue-strict}, if the action of $G$ on $Kv$ is faithful (so in particular, residue-strict valued $G$-fields are strict $G$-fields).
\end{definition}

\begin{lemma}\label{lemma: unique}
If $\Kvs$ is a valued $G$-field, then the valuation $v|_{C_K}$ on $C_K$ extends uniquely to a valuation on $K$.    
\end{lemma}
\begin{proof}
Since the extension $C_K\subseteq K$ is Galois, any two extensions of $v|_{C_K}$ to $K$ are $\operatorname{Gal}\left( K/C_K \right)$-conjugate (see e.g. \cite[Theorem 3.2.14.]{engler}), in particular any such extension is conjugate to $v$. But since $\Kvs$ is a valued $G$-field, the valuation $v$ is invariant under the action of $G\cong \operatorname{Gal}\left(K/C_K\right)$, hence $v$ is the unique extension of $v|_{C_K}$ to $K$.
\end{proof}

\begin{lemma}\label{lemma: ext to residue action faithful}
Every valued $G$-field admits a valued $G$-field extension which is residue-strict.
\end{lemma}
\begin{proof}
Let $\Kvs$ be a valued $G$-field. Consider the function field $L = K\left( X_g\colon g\in G\right)$ with its natural $G$-action, that is $G$ acts on $K$ via $\s$ and $gX_h=X_{gh}$ for $g,h\in G$ . By basic valuation theory (see e. g. \cite[Corollary 2.2.2]{engler}), there is a unique extension $w$ of $v$ to $L$ such that for all $g\in G$:
\begin{enumerate}
    \item $w\left( X_g \right) = 0$,
    \item $\operatorname{res}\left( X_g \right)$ is transcendental over $Kv$.
\end{enumerate}
It is easy to check that the action of $G$ on $\left( L, w \right)$ turns $L$ into a valued $G$-field extension of $K$ and that the action of $G$ on $Lw$ is faithful.
\end{proof}

\begin{lemma}\label{lemma: vK=vC, k^G=c}
Let $K$ be a residue-strict valued $G$-field. Then the extension $C_K\subseteq K$ is defectless, $C_{Kv} = C_Kv$ and $vC_K = vK$.
\end{lemma}
\begin{proof}
By the fundamental inequality for extensions of valued fields we have
$$\left[ K : C_K \right] \ge \left[ vK: vC_K \right] \cdot [ Kv: C_Kv ].$$
Certainly $C_Kv \subseteq C_{Kv}$, thus $\left[ Kv: C_Kv \right] \ge [ Kv: C_{Kv} ]$. Since $G$ acts faithfully on $Kv$ it does so also on $K$, thus $\left[ K: C_K \right] = \left[ Kv: C_{Kv} \right] = | G |
$. Therefore
$$| G | \ge \left[ vK: vC_K \right] \cdot  \left[ Kv: C_Kv \right]  = \left( vK: vC_K \right) \cdot | G | \ge | G|,$$
hence all the inequalities above are equalities. In particular
$$\left[ Kv: C_Kv \right]  = \left[ Kv: C_{Kv} \right]\quad\mbox{ and }\quad\left[ vK: vC_K \right] = 1, $$
thus $C_{Kv} = C_Kv$, $vC_K = vK$ and the extension $C_K\subseteq K$ is defectless.
\end{proof}

Using Lemma \ref{lemma: vK=vC, k^G=c} we can simplify our diagram for residue-strict valued $G$-fields.

\begin{center}
\begin{tikzpicture}[
    >=Latex,
    node distance=0.7cm and 1.5cm,
    every node/.style={font=\normalsize}
]
% top
\node (A){$\bigl(K,v,\sigma\bigr)$};

% first level
\node (C) [below=of A] {$\bigl(C_K,v\bigr)$};
\node (B) [left=of C] {$\bigl(Kv,\sigma\bigr)$};
\node (D) [right= of C] {$vK=vC_K$};

% side nodes
\node (Z) [right=of A] {$\bigl(K,\sigma\bigr)$};
\node (T) [left=of A] {$\bigl(K,v\bigr)$};

% second level
\node (F) [below = of C] {$C_Kv=C_{Kv}$};

% arrows
\draw[->] (A) -- (B);
\draw[->] (A) -- (C);
\draw[->,shorten >=2.8mm] (A) -- (D);
\draw[->] (A) -- (Z);
\draw[->] (A) -- (T);

\draw[->,shorten >=2.8mm] (B) -- (F);
\draw[->] (C) -- (F);
\draw[->] (C) -- (D);

\end{tikzpicture}
\end{center}

The results below is a valuative counterpart of Lemma \ref{lemma: G compositum}.
\begin{lemma}\label{lemma: GVF compositum}
Let $K$ be a valued $G$-field and let $C'$ be a valued field extension of $C_K$, which is linearly disjoint from $K$ over $C_K$. Assume that the valuation on $C'$ extends uniquely to a valuation on $KC'$. Then the action of $G$ has a unique extension to $KC'$ for which the field of constants is $C'$. With this action $KC'$ becomes a valued $G$-field extension of $K$.
\end{lemma}
\begin{proof}
By Lemma \ref{lemma: G compositum} we get a unique $G$-field structure on $KC'$. Since the valuation on $C'$ extends uniquely to a valuation $w$ on $KC'$, this valuation $w$ has to be $\operatorname{Gal}\left( KC'/C'\right)$-invariant. But, by the linear disjointness  $\operatorname{Gal}\left( KC'/C'\right)\cong \operatorname{Gal}\left( K/C_K\right) \cong G$. Therefore $KC'$ with the induced action of $G$ is a valued $G$-field extension of $K$ with the desired properties.
\end{proof}

\begin{lemma}\label{lemma: henselian}
Every valued $G$-field $K$ admits a valued $G$-field extension $L$ for which both $C_L$ and $L$ are henselian.
\end{lemma}
\begin{proof}
Since an algebraic extension of a henselian valued field is henselian we only need to worry about the requirement ``$C_L$ is henselian''. Let $C_K^h$ be the henselianization of $C_K$. The extension $C_K\subseteq C_K^h$ is immediate and by Lemma \ref{lemma: vK=vC, k^G=c} the extension $C_K\subseteq K$ is defectless, thus by Fact \ref{fact: immediate and defectless are ld} we have that $K$ and $C_K^h$ are linearly disjoint over $C_K$. Hence by Lemma \ref{lemma: GVF compositum} the compositum $KC_K^h$ has a structure of a valued $G$-field extension of $K$ with fields of constants $C_K^h$. Hence $K'=KC_K^h$ has all the desired properties.
\end{proof}

We wish to obtain an analogue of Lemma \ref{lemma: ec red to c pure} in the presence of a valuation. To do so, we first have to address the definability of the valuation ring over the valued field of constants\footnote{In the last section we will extend on this topic and discuss the definability of the valution ring over the field of constants as a pure field.}. Recall that we treat any $G$-field $K$ as a (quantifier-free) definable object in $C_K$ (see Remark \ref{rem: K definable in C} and the paragraph below it).

\begin{lemma}\label{cO}
Let $K$ be a strict valued $G$-field. Then there exists a quantifier-free $\lngval$-formula without parameters defining the valuation ring $\cO_{K}$  in $C_K$.
\end{lemma}
\begin{proof}
For simplicity write $C:= C_K$. Let $\overline{\cO}$ denote the integral closure of $\cO_{C}$ in $K$. We aim to show that $\overline{\cO}= \cO_K$. Since valuation rings are integrally closed (see \cite[Theorem 3.1.3]{engler}) we obtain $\overline{\cO}\subseteq \cO_K$. For the reverse inclusion note first that $\left( \cO_{K}\right)^G = \cO_{C}$ since $G$ acts isometrically on $K$. Take any $a\in \cO_{K}$ and consider the following polynomial:
\[
P_a\left( X\right):=\prod_{g\in G}(X-\s_g(a)).
\]
It is a monic polynomial with coefficients from $\cO_{C}$, and  $a$ is one of its roots. Hence $\cO_{K}\subseteq\overline{\cO}$ and thus $\cO_{K}=\overline{\cO}$. Therefore $a\in \cO_{K}$ exactly when $P_a\left( X\right)\in \cO_{C}[X]$. Set $n=|G|=\left[ K:C \right]$ and for $k=1,\dotsc, n$ denote by $s_k\left(x_1,\dotsc , x_n\right)$ the $k$th elementary symmetric polynomial in $n$ variables. We have
\[ P_a\left( X \right) = X^n +\sum_{k=1}^n (-1)^n s_k\left( \s(a)\right)\]
where $\s(a) = \left\langle \s_g(a)\colon g\in G \right\rangle$. By the preceding argument we have $a\in \cO_{K}$ if and only if $s_k\left( \s(a)\right)\in \cO_{C}$, which yields a definition of $\cO_{K}$ in $(C,v)$.
\end{proof}

\begin{lemma}\label{lemma: ec red to c}
Let $K\subseteq L$ be an extension of strict valued $G$-fields. Then $K\ecgval L$ if and only if $C_K\ecval C_L$.
\end{lemma}
\begin{proof}
The ``only if'' part is straightforward, so let us focus on the reverse direction. For the ``if'' part we will argue as in the proof of Lemma \ref{lemma: ec red to c pure}. We already know that $(K,\s)$  is quantifier-free definable in $C_K$ (see Remark \ref{rem: K definable in C}). By the Lemma \ref{cO} we obtain that the valuation on $K$ is quantifier-free definable in $(C,v)$. Hence, we have, that $(K,v,\s)$ is quantifier-free definable in $(C,v)$, which concludes the proof.
\end{proof}

Combining all the results obtained so far we arrive at the following Ax-Kochen/Ershov style result.

\begin{proposition}\label{prop: ec ake for gvf}
Assume that $K$ is residue-strict valued $G$-field and $C_K$ is tame valued field. Then, the following conditions are equivalent:
\begin{enumerate}[itemsep=0.3cm]
    \item $K\ecgval L$,
    \item $C_K\ecval C_L$,
    \item $vC_K\ecoag vC_L$ and $C_Kv\ecrng C_Lv$,
    \item $vC_K\ecoag vC_L$ and $Kv\ecg L$.
\end{enumerate}
\end{proposition}
\begin{proof}
We have $(1)\Longleftrightarrow (2)$ by Lemma \ref{lemma: ec red to c}, $(2)\Longleftrightarrow (3)$ by Fact \ref{fact: ake for tame} and $(3)\Longleftrightarrow (4)$ by Lemma \ref{lemma: ec red to c pure}. 
\end{proof}

The following lemma can be viewed as an analogue of a corollary of Theorem 2.14 from \cite{badplaces}, which provides a construction of an extension of a valued field with prescribed extensions of the value group and the residue field. In the lemma below, we address the analogous problem of constructing a valued $G$-field extension with a prescribed $G$-field residue field extension.

\begin{lemma}\label{lemma: lift G-ext}
Assume that $K$ is residue-strict valued $G$-field. For any $G$-field extension $Kv\subseteq k$ there exists a valued $G$-field extension $K\subseteq K'$ such that $k\cong K'v$ as $G$-field extensions of $Kv$.  
\end{lemma}
\begin{proof}
By \cite[Theorem 2.14]{badplaces} there is a valued field extension $C'$ of $C_K$ for which $C'v\cong k^G$ over $C_Kv$, so we may as well assume that $k^G=C'v$. We aim to show that $K':=KC'$ satisfies the conclusion of the lemma. For readability, we include the following diagrams of extensions appearing in the proof:

\begin{center}
\begin{tikzcd}[every label/.append style={font=\normalsize}]
             & KC'                     &               &                                 &  &    &    &               & k                       &                \\
K \arrow[ru] &                         & C' \arrow[lu] & {} \arrow[rrr,rightarrow,"\operatorname{res}" description, dashed] &  & {} & {} & Kv \arrow[ru] &                          & k^G\arrow[lu] \\
             & C_K \arrow[lu] \arrow[ru] &               &                                 &  &    &    &               & C_Kv \arrow[lu] \arrow[ru] &               
\end{tikzcd}
\end{center}

By Lemma \ref{lemma: vK=vC, k^G=c}
\[
C'v\cap Kv = k^G\cap Kv = (Kv)^G = (Cv)_K=C_Kv.
\]
The extension $C_Kv\subseteq Kv$ is Galois so $C'v\cap Kv =C_Kv $ implies that $Kv$ and $C'v$ are linearly disjoint over $C_Kv$. Therefore:
\begin{enumerate}
    \item the extension $C_K\subseteq K$ is finite and defectless (Lemma \ref{lemma: vK=vC, k^G=c}),
    \item $C'v$ and $Kv$ are linearly disjoint over $C_Kv$,
    \item $vK\cap vC'=vC_K\cap vC'=vC_K$ (Lemma \ref{lemma: vK=vC, k^G=c}),
\end{enumerate}
so by Lemma \ref{lemma: val indep are lin indep} we obtain that $K$ is linearly disjoint from $C'$ over $C_K$, the valuation $v|_K$ extends to $KC'$ uniquely and $Kv.C'v = (KC')v$. By Lemma \ref{lemma: GVF compositum} we get that $KC'$ is (naturally) a valued $G$-field extension of $K$ whose field of constants is $C'$.

Since the action of $G$ on $Kv$ is faithful, so is the action of $G$ on $k$, hence 
\[
\left[ k:C'v\right]=\left[ k:k^G\right] = |G|.
\]
By linear disjointness $\left[ Kv.C'v: C'v\right]=[Kv:C_Kv] = |G|$ therefore
\[
|G|=\left[ k: C'v \right] = \left[ k: Kv.C'v \right] \cdot\left[ Kv.C'v: C'v \right] = \left[ k: Kv.C'v \right]\cdot |G|
\]
hence $\left[ k: Kv.C'v \right]=1$ i.e. $Kv.C'v = k$.
But $Kv.C'v = (KC')v$ hence $(KC')v = k$ as desired.
\end{proof}

\subsection{Existentially closed valued \texorpdfstring{$G$}{G}-fields} Throughout this subsection \textbf{we assume that $\left( K, v,\s\right)$ is an existentially closed valued $G$-field.}

\begin{lemma}\label{lemma: residue action faithful}
$K$ is residue-strict.
\end{lemma}
\begin{proof}
For a fixed element $g\in G$ the property ${\s}_g\neq \operatorname{id}_{Kv}$ is expressible by an existential $\lng$-sentence. Since $K$ is existentially closed as a valued $G$-field, it is enough to construct a valued $G$-field extension of $K$ on which the residue action is faithful, which was done in Lemma \ref{lemma: ext to residue action faithful}.
\end{proof}
Immediately from Lemma \ref{lemma: vK=vC, k^G=c} and Lemma \ref{lemma: residue action faithful} we get the following.
\begin{corollary}\label{cor: vK=vC, k^G=c}
The extension $C_K\subseteq K$ is defectless, $C_{Kv} = C_Kv$ and $vC_K = vK$. 
\end{corollary}

\begin{lemma}\label{lemma: qth powers}
Let $q$ be a prime number. Then every element of $C_K$ (resp. $C_Kv$) is has a $q$th root in $K$ (resp. $Kv$).
\end{lemma}
\begin{proof}
It is enough to prove the statement for $C_K$ and $K$. Assume $a\in C_K$ is not a $q$th power in $K$ and let $b$ be any $q$th root of $a$ inside $\Omega$. It is known (and easy to check) that the polynomial $X^q-a$ is irreducible over $C_K$, thus the extension $C_K\subseteq C_K\left( b \right)$ has degree $q$. Consider the field $L= C_K\left( b\right)\cap K$. Since $C_K\subseteq L\subseteq C_K\left( b \right)$ and $\left[ C_K\left( b \right) : C_K \right] = q$ is prime we have that either $L=C_K$ or $L=C_K\left(b\right)$. In the latter case $b\in L\subseteq K$, so $a$ is a $q$th power in $K$, contrary to the assumption. Thus $L= C_K$ and since $C_K\subseteq K$ is a Galois extension we have that $K$ and $C_K\left( b \right)$ are linearly disjoint over $C_K$. Thus by Lemma \ref{lemma: GVF compositum} the compositum $KC_K(b)=K(b)$ becomes a valued $G$-field extension of $K$. But $K\subseteq K(b)$ is a proper algebraic extension, which contradicts the fact that $K$ is existentially closed as a valued $G$-field.
\end{proof}

Immediately from Lemma \ref{lemma: perfect constants} and Lemma \ref{lemma: qth powers} we get the following.
\begin{corollary}\label{cor: perfect}
The fields $K, C_K, Kv, C_Kv$ are perfect.
\end{corollary}

\begin{lemma}\label{lemma: divisible}
The group $vC_K=vK$ is divisible.
\end{lemma}
\begin{proof}
Pick any $a\in C_K$ and any prime number $q$. By Lemma \ref{lemma: qth powers} we have that there is some $b\in K$ such that $b^q = a$ and hence $v\left( a \right) =qv\left( b \right)\in q\cdot vK$. On the other hand, by Corollary \ref{cor: vK=vC, k^G=c} we have $vK = vC_K$, so $v\left( a \right)$ is $q$-divisible in $vC_K$. Hence $vC_K$ is a divisible group.
\end{proof}

\begin{lemma}\label{lemma: tame}
The valued fields $C_K$ and $K$ are tame.
\end{lemma}
\begin{proof}
By \cite[Lemma 2.13]{tame} an algebraic extension of a tame fields is again a tame field, thus it is enough to show that $C_K$ is tame.
By Corollary \ref{cor: perfect} and Lemma \ref{lemma: divisible} the only thing left to prove is that $C_K$ is algebraically maximal. Let $C_K\subseteq C'$ be an algebraic immediate extension. We aim to show that $C'=C_K$. By Fact \ref{fact: immediate and defectless are ld} and Corollary \ref{cor: vK=vC, k^G=c} we have that $K$ and $C'$ are linearly disjoint over $C_K$, so by Lemma \ref{lemma: GVF compositum} and Fact \ref{fact: immediate and defectless are ld} again the field $KC'$ becomes naturally a valued $G$-field extension of $K$. Since $K$ is existentially closed as a valued $G$-field and $KC'$ is algebraic over $K$ we have that $KC'=K$ and thus $C'=C_K$ as desired.
\end{proof}

\begin{lemma}\label{lemma: kv ec}
The residue field $Kv$ is existentially closed as a $G$-field.
\end{lemma}
\begin{proof}
Let $Kv\subseteq k$ be a $G$-field extension. Pick $K'$ as in Lemma \ref{lemma: lift G-ext}. Since $K$ is existentially closed in $K'$ as a valued $G$-field it follows that $Kv$ is existentially in $K'v\cong k$ as a $G$-field. Thus $Kv$ is existentially closed as a $G$-field.
\end{proof}

In the following corollary we summarize some of the results obtained so far and list them in a way that makes following the proof of upcoming Theorem \ref{thm: main ake} more pleasant. 

\begin{corollary}\label{cor: ec satisfy gtcvf}
Assume that $K$ is an existentially closed valued $G$-field with field of constants $C_K$. Then:
\begin{enumerate}
\item $C_K$ is a tame valued field,
\item $vK$ is existentially closed as an ordered abelian group,
\item $Kv$ is existentially closed as a $G$-field.
\end{enumerate}
\end{corollary}
\begin{proof}
By Lemma \ref{lemma: tame} the valued field $C_K$ is tame, by Lemma \ref{lemma: divisible} the value group $vK$ is divisible, hence existentially closed as an ordered abelian group, and finally by Lemma \ref{lemma: kv ec} the residue field $Kv$ is existentially closed as a $G$-field.
\end{proof}

\begin{proposition}
The field $K$ is not PAC. In particular, the $G$-field $K$ is not existentially closed.
\end{proposition}
\begin{proof}
Assume $K$ is PAC. By the Frey–Prestel theorem (see \cite[Corollary 11.5.5]{fields}) we have that the henselianization of the valued field $K$ is the separable closure of $K$ and that the residue field $Kv$ is separably closed. But $G$ acts faithfully on $Kv$, so by the Artin-Schreirer theorem we have that:
\begin{enumerate}
    \item $G = C_2$ is the cyclic group of order two,
    \item $Kv$ is algebraically closed of characteristic zero,
    \item the fixed field $C_Kv$ is real closed.
\end{enumerate}
By \cite[Example 3.5]{GTCF} an existentially closed $C_2$-field cannot be algebraically closed. The ``In particular'' part follows, since existentially closed $G$-fields are PAC by \cite[Theorem 3.29]{GTCF}.
\end{proof}

Whilst investigating the model-theoretic properties of existentially closed fields in Section \ref{NTP2} the following notion will be of great importance.

\begin{definition}[Introduction of \cite{BlaszczokKuhlmann2017}]
Let $(K,v)$ be a valued field of residue characteristic $p$. We say that $(K,v)$ is \textbf{Kaplansky} if either
\begin{enumerate}
\item $p=0$ and $vK$ is divisible, or
\item $p>0$, $vK$ is $p$-divisible and $Kv$ has no finite extension of degree divisible by $p$.
\end{enumerate}
\end{definition}

\begin{lemma}\label{lemma: C is benign}
Assume $\operatorname{char} Kv\nmid |G|$. Then the valued field $(C_K,v)$ is Kaplansky.
\end{lemma}
\begin{proof}
Let $p:=\operatorname{char} Kv$. If $p=0$ then we are done by Lemma \ref{lemma: divisible}, so assume that $p>0$. Let $F$ be a finite algebraic extension of $C_Kv$ of degree $d$ divisible by $p$. The field $C_K v$ is perfect by Lemma \ref{lemma: perfect constants} therefore $F$ is a separable extension of $C_Kv$ and by passing to the normal closure we may assume that $C_Kv\subseteq F$ is a Galois extension. Let $C_Kv\subseteq F'\subseteq F$ be a subextension of $F$ such that $[F':C_Kv]=p^k$ for some integer $k$ and $[F:F']=d$ is coprime to $p$.\footnote{Simply take a Sylow $p$-subgroup $H\le\operatorname{Gal}(F/C_Kv)$ and set $F':=F^H$.} Then $F'$ is linearly disjoint from $Kv$ over $C_Kv$ so by Lemma \ref{lemma: G compositum} the compositum $F'Kv$ is naturally a finite $G$-field extension of $Kv$. But $Kv$ is existentially closed as a $G$-field by Lemma \ref{lemma: kv ec} which forces $k=0$.
\end{proof}

\begin{remark}
    Note that assumption $\operatorname{char} Kv\nmid |G|$ is crucial for Lemma \ref{lemma: C is benign} since the extension $C_Kv\subseteq Kv$ is of degree $| G|$. 
\end{remark}

\section{Model companion}

We are ready to prove the main result of this paper.
\begin{theorem}\label{thm: main ake}
Assume $G$ is a finite nontrivial group and let $K$ be a valued $G$-field. Then $K$ is existentially closed as a valued $G$-field if and only if
\begin{enumerate}
    \item $C_K$ is tame,
    \item $vK$ is existentially closed as an ordered abelian group,
    \item $Kv$ is existentially closed as a $G$-field.
\end{enumerate}
\end{theorem}
\begin{proof}
The ``only if'' direction is exactly Corollary \ref{cor: ec satisfy gtcvf}. For the other direction let $K$ be a valued $G$-field such that $C_K$ is tame, $vK$ is divisible and $Kv$ is existentially closed as a $G$-field (in particular, $K$ is residue-strict). Let $K'$ be the valued $G$-field extension of $K$. By Proposition \ref{prop: ec ake for gvf} in order to show that $K\ecgval K'$ it is enough to show that
\begin{enumerate}
    \item $vC\ecoag vC_{K'}$ and
    \item $Kv \ecg K'v$.
\end{enumerate}
But this follows since both $vC_K=vK$ and $(Kv,\s)$ are existentially closed.
\end{proof}

By Subsection 7.1 of \cite{tame} we have that the class of tame fields is elementary. Let $\gtcvf$ be the $\lng$-theory expressing the following properties:
\begin{enumerate}
    \item $C_K$ is tame,
    \item $\left( Kv, \s\right)$ is a model of $G\operatorname{TCF}$.
    \item $vK$ is divisible.
\end{enumerate}
Using this definition and the fact that the $\lng$-theory of valued $G$-fields is inductive, Theorem \ref{thm: main ake} translates into the following.
\begin{theorem}\label{thm: mc}
The theory $\gtcvf$ is the model companion of the theory of valued $G$-fields in the language $\lng$.
\end{theorem}

 Now, we wish to provide an example of a model of $\gtcvf$. 

\begin{example}
     In Example 3.11 \cite{GTCF} Kowalski and Hoffmann constructed a model of $(\Z/n\Z)$TCF. Let us denote this model by $K$ and its field of constants by $C$. Let $\Gamma$ be a divisible ordered abelian group. Consider the Hahn series field $C((t^\Gamma))$ with a standard valuation (with the value group $\Gamma$ and the residue field $C$). Recall, that  by the classical result of Krull from \cite{Krull}\footnote{To be precise, it was generalized by Poonen in \cite{Poonen} to apply to the broader context, like the one in our example.} the valued field $C((t^\Gamma))$ is algebraically maximal. By Remark 3.20 in \cite{GTCF} the residue field $C$ is perfect. Finally, by our choice of divisible $\Gamma$, the valued field $C((t^\Gamma))$ is tame.
     
    Note that since we assumed $\Gamma$ to be divisible, so in particular torsion-free, $C$ is relatively algebraically closed in $C((t^\Gamma))$. Therefore we obtain that $K$ and $C((t^\Gamma))$ are linearly disjoint over $C$, and by the usual argument we can extend action of $\Z/n\Z$ from $K$ to compositum field $C((t^\Gamma))K=K((t^\Gamma))$. Since the valuation from $C((t^\Gamma))$ to $K((t^\Gamma))$ extends uniquely, we obtain on $K((t^\Gamma))$ the structure of valued $\Z/n\Z$-field with the following properties:
    \begin{enumerate}
        \item the field of constants $C((t^\Gamma))$ is tame,
        \item the residue field $K$ is model of $\Z/n\Z$TCF,
        \item the value group $\Gamma$ is divisible,
    \end{enumerate}
    so, by Theorem \ref{thm: main ake} $K((t^\Gamma))$ is a model of $(\Z/n\Z)$TCVF. 

\end{example}

\section{Neostability of \texorpdfstring{$\gtcvf$}{GTCVF}}\label{NTP2}
Let us fix $p$ which is either zero or a prime number and denote by $\gtcvf_p$ the theory whose models are $K\models\gtcvf$ of residue characteristic $p$. In this section we will show that (provided that $p$ does not divide the order of $G$) the theory $\gtcvf_p$ is of finite burden, so in particular it has $\ntp$. The proof is a rather straightforward application of: 
\begin{enumerate}
    \item an Ax-Kochen/Ershov style result for calculating burden of valued fields (proven by Touchard in \cite{Pierre}), and
    \item the supersimplicity of $\gtcf$ (proven by Hoffmann and Kowalski in \cite{GTCF}).
\end{enumerate}
We refer to \cite{ntp2} for the necessary background on $\ntp$ theories. We denote the burden of a definable set $X$ by $\bdn X$.

\textbf{Through this section we will assume that $p$ does not divide the order of the group $G$.}

\begin{remark}
There is no hope for $\gtcvf_p$ to be simple or have NIP. Take any valued $G$-field $K\models \gtcvf_p$. If $\gtcvf_p$ were simple (resp. have NIP) then $vK$ (resp. $Kv$) would be both simple and have NIP, hence stable, which it is not.
\end{remark}

\begin{fact}[Theorem 3.12, Corollary 3.13 in \cite{Pierre}]\label{bdn}
Let $(K,v)$ be an algebraically maximal Kaplansky valued field with perfect residue field. Then 
\[
\bdn(K)=\max\{\bdn(Kv),\bdn(vK)\}.
\]
\end{fact}

\begin{theorem}\label{thm: ntp2}
The (theory of the) valued field $C_K$ is inp-minimal i.e. it has burden equal one.
\end{theorem}
\begin{proof}
The valued field $C_K$ is tame valued field, thus Lemma \ref{lemma: C is benign} implies that $C_K$ satisfies the assumptions of Fact \ref{bdn}. Since $vC_K$ is a divisible ordered abelian group we have that $\bdn(vC_K)=1$. On the other hand, the theory of $C_Kv$ is supersimple and of SU-rank one (see \cite[Proposition 4.8]{GTCF}). Therefore \cite[Theorem 4.4.8]{kim} yields that $\bdn(C_Kv)=1$. Thus, by Fact \ref{bdn}, the burden of $C_K$ equals one, as desired.
\end{proof}

\begin{corollary}
The theory $\gtcvf_p$ is of finite burden with burden bounded by a number depending only on $| G|$.
\end{corollary}
\begin{proof}
Take any valued $G$-field $K\models \gtcvf_p$. The structure $\Kvs$ is interpretable in its valued field of constants $(C_K,v)$ (see Remark \ref{rem: K definable in C}). Since burden is sub-multiplicative (see \cite[Corollary 2.6]{ntp2}) the burden of $(K,v,\s)$ is finite by Theorem \ref{thm: ntp2} and bounded in terms of $| G|$.
\end{proof}

\begin{remark}
We do not know what happens if $p$ divides the order of $G$. We suspect that $\gtcvf_p$ is still $\ntp$, but possibly of infinite burden. It might be an interesting test case for extending the results of \cite{Pierre} (e.g. Fact \ref{bdn} in its general form) to some natural class of valued fields including models of $\gtcvf_p$.
\end{remark}

\section{Definability of the valuation ring}
Fix a model $K\models\gtcvf$ and set $C:=C_K$. In this section we prove that $\cO_C$ is existentially definable in the pure field $C$ by applying results of Fehm from \cite{fehm}.

Let us recall the set-up from \cite{fehm}. For a polynomial $f\in C\left[ X\right]$ consider the following formulas
\begin{enumerate}
    \item $\varphi_f(x) \equiv (\exists y,z,y_1,z_1)( x=y_1-z_1 \wedge y_1f(y)=1 \wedge z_1f(z)=1 )$,
    \item $\psi_f(x) \equiv (\exists y,z,y_1,z_1)(x=0\vee(x=y_1z_1\wedge y_1f(y)=1\wedge z_1f(z)=1))$,
    \item $\eta_f(x)\equiv(\exists u,t)(x=u+t\wedge\varphi_f(u)\wedge\psi_f(t))$,
\end{enumerate}
so that (provided that $f$ has no root in $C$)
\begin{enumerate}
    \item $\varphi_f(C) = f(C)^{-1}-f(C)^{-1}$,
    \item $\psi_f(C) = f(C)^{-1}\cdot f(C)^{-1}$,
    \item $\eta_f(C) = \varphi_f(C) + \psi_f(C)$.
\end{enumerate}
With this notation we have the following.

\begin{fact}[Proposition 3.3. in \cite{fehm}]\label{fact: fehm}
Let $C$ be a valued field and let $f\in \cO_C[ X]$ be a monic polynomial such that $\overline{f}\in Cv[X]$ is (absolutely) square-free and has no zero in $Cv$. Assume that $Cv$ is PAC and there is $a\in \cO_C$ such that $f'(a)\not\in \mathfrak{m}$. Then $\eta_f\left( C \right) = \cO_C$.
\end{fact}

Let $b\in Kv$ be such that $(Cv)[b]=Kv$, let $f_0$ be the minimal polynomial of $b$ over $Cv$ and let $f\in \cO_C[X]$ be a monic lift of $f_0$. Finally, let $\overline{c}$ be the tuple of coefficients of $f$. With this notation we have the following.

\begin{proposition}\label{prop: def val}
Let $\left( K, \s, v \right)$ be a model of $\gtcvf$ and denote by $C$ the fields of constants of $K$. Then $\cO_C$ is existentially definable in the pure field $C$. Moreover, there is an existential formula $\eta(x)$ over the empty set  in the language $\mathcal{L}_{\operatorname{rng}}(\overline{c})$ such that for any $K'\models \gtcvf$ extending $K$ we have $\eta( C' ) =\cO_{C'}$ where $C'$ is the field of constants of $K'$.
\end{proposition}
\begin{proof}
Since $Cv$ is perfect by Corollary \ref{cor: perfect} the polynomial $f_0$ is square-free and $f'_0$ is not identically zero. By Lemma \ref{lemma: kv ec} the field $Kv$ is a model of $\gtcf$, so by Fact \ref{gtcf} the field $Cv$ is PAC so in particular infinite. Therefore there is some $a\in \cO_C$ such that $\operatorname{res}(a)$ is not a root of $f_0$, i.e. $f'(a)\not\in \mathfrak{m}$ and thus $\eta_f\left( C \right) = \cO_C$ by Fact \ref{fact: fehm}. For the ``Moreover...'' part it is enough to note that for any $K'\models \gtcvf$ extending $K$ the polynomial $f$ still satisfies the assumptions of Fact \ref{fact: fehm} since $K$ is existentially closed in $K'$ by Theorem \ref{thm: mc}.
\end{proof}
\begin{remark}
Due to Remark \ref{rem: K definable in C} and Proposition \ref{prop: def val}, the structure of $K$ as a valued $G$-field (i.e. the valuation and the action of $G$) is interpretable in the pure field $C_K$ alone.
\end{remark}

\begin{proposition}\label{prop: mc field}
The $\lngring(\bar{c})$-theory of $C$ is model complete.
\end{proposition}
\begin{proof}
Let $C\subseteq C'$ be a pure field extension with $C'\equiv C$ as $\lngring(\bar{c})$-structures. We need to show that $C\prec C'$.
\begin{claim}
The fields $K$ and $C'$ are linearly disjoint over $C$.
\end{claim}
\begin{claimproof}
Since $C'\equiv C$ as $\lngring(\bar{c})$-structures we have $K\cap C'= C$ and since $C\subseteq K$ is a Galois extension the claim follows.
\end{claimproof}
By Claim we get that $\operatorname{Gal}(KC'/C') \cong \operatorname{Gal}(K/C)$ with the isomorphism being induced by the restriction map, hence $K':=KC'$ becomes a $G$-field extension of $K$ with $C_{K'}=C'$.

Take $\eta(x)$ as in Proposition \ref{prop: def val}. The information
\begin{center}
$\eta(C)$ is a valuation ring of $C$
\end{center}
is part of the $\lngring(\bar{c})$-theory of $C$ and since $C'\equiv C$ as $\lngring(\bar{c})$-structures we have
\begin{center}
$\eta(\cO')$ is a valuation ring of $C$.
\end{center}
Therefore $K'$ is a valued field with an action of $G$. The $G$-invariance of $\cO$ is also a part of the $\lngring(\bar{c})$-theory of $C$, hence $K'$ is a valued $G$-field. Since $C\equiv C'$ we have $K\equiv K'$ as $G$-fields. Since $K\models \gtcvf$ and $\gtcvf$ is model complete we have $K\ecgval K'$ and hence $C\prec C'$ even as $\lngring(\bar{c})$-structures.
\end{proof}

\begin{example}\label{ex: c in Q}
Let us construct an example of $K\models \gtcvf$ for which we can take $\overline{c}$ belonging to $\mathbb{Q}$. For simplicity we will take $G=C_2$ (the cyclic group of order two) The element $2\in \mathbb{Q}_3$ has no square root in $\mathbb{Q}_3$ so $K_0=\mathbb{Q}_3\left[\sqrt{2}\right]$ is a proper extension of $C_0=\mathbb{Q}_3$. The field $C_0$ is henselian with respect to the $3$-adic valuation $v_3$, thus $v_3$ extends uniquely to a valuation $v$ on $K_0$. Hence $\left( K_0 , v\right)$ with its natural $G$-action is a valued $G$-field and $C_{K_0} = C_0$. The theory of valued $G$-field is inductive so we can extend $K_0$ to a model $K\models \gtcvf$. The element $\sqrt{2}\in K$ is not fixed by $G$, thus the extension $C_K\subseteq C_K\left[ \sqrt{2} \right]$ is proper. Since $\left[ K : C_K \right] = | G| = 2$ we have $K=C_K\left[ \sqrt{2} \right]$. Using the notation introduced right before Proposition \ref{prop: def val} we can take
\begin{enumerate}
    \item $b= \operatorname{res}\sqrt{2}\in Kv$,
    \item $f_0(X) = X^2-2\in (C_Kv)\left[ X\right]$,
    \item $f(X) = X^2-2\in C_K\left[ X\right]$,
\end{enumerate}
hence $\overline{c} = \left( -2, 0, 1\right)\in \mathbb{Z}^3$. In this case the formula $\eta(x)$ in Proposition \ref{prop: def val} can be taken to be an $\lngring$-formula over the empty set. Proposition \ref{prop: mc field} yields that the pure field $C$ is model complete.
\end{example}

\begin{remark}\label{rem: mc fields}
Proposition \ref{prop: mc field} yields interesting new examples of a model complete fields. If $\overline{c}$ is in the prime field of $C$ then $C$ is a pure field which:
\begin{enumerate}
    \item is model complete,
    \item is inp-minimal,
    \item carries a 0-definable henselian valuation,
    \item interprets a supersimple PAC field.
\end{enumerate}
In particular, the field $C$ constructed in Example \ref{ex: c in Q} satisfies all the above properties.
\end{remark}

\bibliographystyle{amsplain}
\bibliography{biblio.bib}

\end{document}